\documentclass[11pt,oneside]{amsart}
\usepackage[T1]{fontenc}
\usepackage[latin9]{inputenc}
\usepackage{geometry}
\usepackage{verbatim}
\usepackage{textcomp}
\usepackage{amsthm}
\usepackage{amstext}
\usepackage{amssymb}
\usepackage{esint}
\usepackage{mathrsfs}
\usepackage{wasysym}
\usepackage{nomencl}
\usepackage[all]{xy}
\usepackage{cancel} % strikeouts
\usepackage{color}

\makeatletter
\numberwithin{equation}{section}
\numberwithin{figure}{section}
\theoremstyle{plain}
\newtheorem{thm}{\protect\theoremname}
  \theoremstyle{definition}
  \newtheorem{defn}[thm]{\protect\definitionname}
  \theoremstyle{plain}
  \newtheorem{prop}[thm]{\protect\propositionname}
  \theoremstyle{plain}
  \newtheorem{lem}[thm]{\protect\lemmaname}
  \theoremstyle{remark}
  
  \theoremstyle{definition}
  \newtheorem{example}[thm]{\protect\examplename}
   \theoremstyle{definition}
  \newtheorem{cor}[thm]{\protect\corollaryname}
\theoremstyle{definition}
  \newtheorem{defp}[thm]{\protect\definitionpropositionname}
\theoremstyle{definition}
  \newtheorem{defc}[thm]{\protect\definitioncorollaryname}
\def\R{{\mathbb R}}

\def\bQ{{\mathbb Q}}

\def\bN{{\mathbb N}}

\def\SH{{\mathbf{Sh}}}
\def\QSH{{\mathbf{QSh}}}
\def\QSUL{{\mathbf{QSul}}}

\def\QSym{{\mathbf {QSym}}}

\def\N{{\mathbf N}}

\def\bQ{{\mathbf Q}}

\def\Sh{{\mathbf {Sh}}}

\def\shuff#1#2{\mathbin{
      \hbox{\vbox{\hbox{\vrule \hskip#2 \vrule height#1 width 0pt}\hrule}\vbox{\hbox{\vrule \hskip#2 \vrule height#1 width 0pt\vrule }\hrule}}}}
\def\shuffl{{\mathchoice{\shuff{5pt}{3.5pt}}{\shuff{5pt}{3.5pt}}{\shuff{3pt}{2.6pt}}{\shuff{3pt}{2.6pt}}}}
\def\shuffle{\, \shuffl \,}
\def\dshuffle{\overline\shuffle}

\def\qshuffl{{\mathchoice{\shuff{5pt}{3.5pt}\hspace{-2.9mm}-}{\shuff{5pt}{3.5pt}\hspace{-2.9mm}-}
{\shuff{3pt}{2.6pt}\hspace{-2.2mm}-}{\shuff{3pt}{2.6pt}\hspace{-2.2mm}-}}}
\def\qshuffle{\,\qshuffl\,}
\makeatother

\usepackage[english]{babel}
  \providecommand{\definitionname}{Definition}
  \providecommand{\examplename}{Example}
  \providecommand{\lemmaname}{Lemma}
  \providecommand{\propositionname}{Proposition}
  \providecommand{\remarkname}{Remark}
\providecommand{\theoremname}{Theorem}

\providecommand{\corollaryname}{Corollary}
\providecommand{\definitionpropositionname}{Definition-Proposition}
\providecommand{\definitioncorollaryname}{Definition-Corollary}
\author{Claudia Malvenuto}
\address{Dipartimento di Matematica\\ Sapienza Universit\`a
di Roma\\ P.le A. Moro 2\\ 00185, Roma, Italy\\
email {claudia@mat.uniroma1.it}}

\author{Fr\'ed\'eric Patras}
\address{UMR 7351 CNRS\\
        		Universit\'e de Nice\\
        		Parc Valrose\\
        		06108 Nice Cedex 02
        		France\\
        		email {patras@unice.fr}}
\begin{document}

%%%%%%%%%%%%%%%%%%%%%%%%%%%%%%%%%%%%%%%%%
%%%%%%%%%%%%%%%%%%%%%%%%%%%%%%%%%%%%%%%%%

\selectlanguage{english}

\title[Resummation of MZVs]{Symmetril moulds, generic group schemes, resummation of MZVs}

\date{}

\maketitle
\begin{abstract} The present article deals with various generating series and group schemes (not necessarily affine ones) associated with MZVs. Our developments are motivated by Ecalle's mould calculus approach to the latter. We propose in particular a Hopf algebra--type encoding of symmetril moulds and introduce a new resummation process for MZVs.
\end{abstract}

\section*{Introduction}\label{sect:intro}
 
Motivated by the study of multiple zeta values (MZVs), Jean Ecalle has introduced various combinatorial notions such as the ones of ``symmetral moulds'', ``symmetrel moulds'', ``symmetril moulds'' or ``symmetrul moulds'' \cite{ARI2003,cresson}. The first two are well-understood classical objects: they are nothing but characters on the shuffle algebra, resp. the quasi-shuffle algebra over the integers, both isomorphic to the algebra $\QSym$ of Quasi-symmetric functions. These two notions are  closely related to the interpretation of properly regularized MZVs as real points of two prounipotent affine group schemes (associated respectively to the integral and power series representations of MZVs), whose interactions through double shuffle relations has given rise to the  modern approach to MZVs (by Zagier, Deligne, Ihara, Racinet, Brown, Furusho and may others) \cite{cartier, furu, ihara, racinet}

Although fairly natural from the point of view of MZVs (the resummation of MZVs info suitable generating series gives rise to a symmetril mould), the notion of symmetrility is more intriguing and harder to account for using classical combinatorial Hopf algebraic tools. 

The aim of this article is accordingly threefold.
We first show that Ecalle's mould calculus can be interpreted globaly, beyond the cases of symmetral and symmetrel moulds, as a rephrasing of the theory of MZVs into the framework of prounipotent groups. However, these not necessarily associated to affine group schemes (that is, to groups whose elements are characters on suitable Hopf algebras), at least in our interpretation and indeed, to account for symmetrility we introduce a new class of functors from commutative algebras to groups refered to as generic group schemes (because the elements of these groups are characters on suitably defined ``generic'' Hopf algebras). Second, we focus on this notion of symmetrility, develop systematic foundations for the notion and prove structure theorems for the corresponding algebraic structures. Third, we interpret Ecalle's resummation of MZVs by means of formal power series as the result of a properly defined Hopf algebra morphism. This construction is reminiscent in many aspects of the resummation of the various Green's functions in the functional calculus approach to quantum field theory or statistical physics, see e.g. \cite{PS}), 
This approach leads us to introduce a new resummation process, different than Ecalle's. The new process is more complex combinatorially but more natural from the group and Lie theoretical point of view: indeed, it encodes MZVs into new generating series that behave according to the usual combinatorics of tensor bialgebras and their dual shuffle bialgebras.

In the process, we introduce various Hopf algebraic structures that, besides being motivated by the mould calculus approach to MZVs, seem to be interesting on their own from a combinatorial algebra point of view.

We refer the readers not acquainted with classical arguments on the theory of MZVs to Cartier's Bourbaki seminar \cite{Cartier} that provides a short and mostly self contained treatment of the key notions. 

%%%%%%%%%%%%%%%%%%%%%%%%%%%%%%%%%%%%%%%
\section*{Acknowledgements} The authors acknowledge support from ICMAT, Madrid, and from the grant CARMA, ANR-12-BS01-0017.

\section{Hopf algebras}\label{sec:ha}

We recall first briefly the definition of a Hopf algebra and related notions. The reader is refered to \cite{Cartier} for details. All the maps we will consider between vector spaces will be assumed to be linear excepted if otherwise stated explicitely.
We will be mostly interested in graded or filtered connected Hopf algebras, and restrict therefore our presentation to that case. 

Let $H=\bigoplus_{n\in \N}H_n$ be a graded vector space over a field $k$ of characteristic zero. We will always assume that the $H_n$ are finite dimensional.
We write $H_{\leq n}:=\bigoplus\limits_{m\leq n}H_m$, $H_{\geq n}:=\bigoplus\limits_{m\geq n}H_m$ and $H^+:=\bigoplus_{n\in \N^\ast}H_n$.
The graded vector space $H$ is said to be connected if $H_0\cong k$. 

An associative and unital product $\mu:H\otimes H\to H$ on $H$ (also written $h\cdot h':=\mu(h\otimes h')$) with unit map $\eta: k\rightarrow H_0\subset H$ (so that for any $h\in H$ and $\eta(1)=:{\mathbf 1}\in H_0$, ${\mathbf 1}\cdot h=h\cdot {\mathbf 1}=h$) makes $H$ a graded (resp. filtered) algebra if, for any integers $n,m$, $\mu(H_n\otimes H_m)\subset H_{n+m}$ (resp. $\mu(H_n\otimes H_m)\subset H_{\leq n+m}$). 

Dualizing, a coassociative and counital coproduct $\Delta:H\to H\otimes H$ on $H$ (also written using the abusive but useful Sweedler notation $\Delta(h)=h^{(1)}\otimes h^{(2)}$) with counit map $\nu:H\rightarrow k$ (with $\nu$ the null map on $H^+$) makes $H$ a graded  coalgebra if, for any integer $n$, $\Delta(H_n)\subset \bigoplus\limits_{p+q=n}H_p\otimes H_q=:(H\otimes H)_n$. The coproduct (resp. the coalgebra $H$) is cocommutative if for any $h\in H$, $h^{(1)}\otimes h^{(2)}=h^{(2)}\otimes h^{(1)}$.

Recall that the category of associative unital algebras is monoidal: the tensor product of two associative unital algebras is a unital associative algebra. Assume that ($\mu$,$\eta$) and ($\Delta$,$\nu$) equip $H$ with the structure of an associative unital algebra and coassociative counital coalgebra: 
they equip $H$ with the structure of a bialgebra if furthermore $\Delta$ and $\nu$ are maps of algebras (or equivalently $\mu$ and $\eta$ are map of coalgebras). The bialgebra $H$ is called a Hopf algebra if furthermore there exists a endomorphism $S$ of $H$ (called the antipode) such that
\begin{equation}
\mu\circ (Id \otimes S) \circ \Delta  = \mu\circ (S \otimes Id) \circ \Delta  = \eta\circ\nu=: \varepsilon .
\end{equation}
A bialgebra or a Hopf algebra is graded (resp. filtered) if it is a graded algebra and coalgebra (resp. a filtered algebra and a graded coalgebra).
Graded and filtered connected bialgebras are automatically equipped with an antipode and are therefore Hopf algebras, and the two notions of Hopf algebras and bialgebras identify in that case, see e.g. \cite{Cartier} for the graded case, the filtered one being similar. This observation will apply to the bialgebras we will consider.

\begin{example} The first example of a bialgebra occuring in the theory of MZVs is $\QSym$, the quasi-shuffle bialgebra over the integers $\N^\ast$. The underlying graded vector space is the vector space over the sequences of integers (written as bracketed words) $[n_1...n_k]$. The bracket notation is assumed to behave multilinearly: for example, for two words $n_1\dots n_k$, $m_1\dots m_l$ and two scalars $\lambda,\beta$
$$[\lambda \ n_1\dots n_k+\beta\  m_1\dots m_l]=\lambda [n_1\dots n_k]+\beta [m_1\dots m_l].$$

The words of lenght $k$ span the degree $k$ component of $\QSym$ (another graduation is obtained by defining the word $[n_1...n_k]$ to be of degree $n_1+\dots +n_k$). The graded coproduct is the deconcatenation coproduct:
$$\Delta([n_1...n_k]):=\sum\limits_{i=0}^k[n_1\dots n_i]\otimes [n_{i+1}\dots n_k].$$
The unital ``quasi-shuffle'' product $\qshuffle$ is the filtered product defined inductively by (the empty word identifies with the unit):
$$[n_1...n_k]\qshuffle [m_1\dots m_l]:= [n_1(n_2...n_k\qshuffle m_1\dots m_l)]+$$
$$[m_1(n_1...n_k\qshuffle m_2\dots m_l)]+
 [ (n_1+m_1)(n_2...n_k\qshuffle m_2\dots m_l)].$$
  For example, $$[35]\qshuffle [1]=[3(5\qshuffle 1)+1(35)+45]=[351]+[315]+[36]+[135]+[45].$$
  Notice that, here and later on, we use in such formulas the shortcut notation $[3(5\qshuffle 1)]$ for the concatenation of $[3]$ with $[5\qshuffle 1]$.

Algebra characters on $\QSym$ (i.e. unital multiplicative maps from $\QSym$ to a commutative unital algebra $A$) are called by Ecalle \it symmetrel moulds\rm . 
The convolution product of linear morphisms from $\QSym$ to $A$,
$$f\ast g:=m_A \circ (f\otimes g)\circ \Delta,$$ where $m_A$ stands for the product in $A$,
equips 
the set $G_\QSym(A)$ of $A$-valued characters with a group structure. Since $\QSym$ is a filtered connected commutative Hopf algebra,
the corresponding functor $G_\QSym$ is (by Cartier's correspondence between group schemes and commutative Hopf algebras over a field of characteristic 0) a prounipotent affine group scheme.
 Properly regularized MZVs are real valued algebra characters on $\QSym$ and probably the most important example of elements in $G_\QSym({\mathbf R})$ \cite{Cartier}.

The quasi-shuffle bialgebra $\QSH (B)$ over an arbitrary commutative algebra $(B,\times)$ is defined similarly: the underlying vector space is $T(B):=\bigoplus\limits_{n\in{\mathbf N}}B^{\otimes n}$, the coproduct is the deconcatenation coproduct and the product is defined recursively by (we use a bracketed word notation for tensor products):
$[b_1\dots b_k]:=b_1\otimes ...\otimes b_k$,
$$[b_1\dots b_k]\qshuffle [c_1 \dots  c_l]:= [b_1(b_2...b_k\qshuffle c_1\dots c_l)]+$$
$$[c_1(b_1...b_k\qshuffle c_2\dots c_l)]+
 [ (b_1\times c_1)(b_2...b_k\qshuffle c_2\dots c_l)].$$

  \end{example}
\begin{example}\label{ex2} The second example arises from the integral representation of MZVs. The corresponding graded vector space $T(x,y)$ is spanned by words in two variables $x$ and $y$. The lenght of a word defines the grading.
The coproduct is again the deconcatenation coproduct acting on words. The product $\shuffle$ is the shuffle product, defined inductively on sequences by
$$a_1...a_k\shuffle b_1\dots b_l:=a_1(a_2...a_k\shuffle b_1\dots b_l)+b_1(a_1...a_k\shuffle b_2\dots b_l).$$

The Hopf algebra $T(x,y)$ is called the shuffle bialgebra over the set $\{x,y\}$.
Properly regularized MZVs are algebra characters on $T(x,y)$ (or on subalgebras thereof), but the regularization process fails to preserve simultaneously the shuffle and quasi-shuffle products \cite{Cartier}.

\it Shuffle bialgebras \rm over arbitrary sets $X$ are defined similarly and denoted $\SH (X)$ (see \cite{Foissy2}). In the mould calculus terminology, a character on $\SH (X)$ is called a \it symmetral mould\rm . The shuffle bialgebra over $\N^\ast$, $\SH(\N^\ast)$, is written simply $\SH$ and will be called the \it integer shuffle bialgebra\rm . It is isomorphic to $\QSym$ as a bialgebra \cite{Hoffman}.
\end{example}

\begin{example} \it Rota-Baxter quasi-shuffle bialgebras\rm . This third example departs from the two previous ones in that it is not a classical one but already illustrates a leading idea of mould calculus, namely: the application of fundamental identities of integral calculus to word-indexed formal power series. We refer e.g. to \cite{eg} and to the survey \cite{epGazette} for an overview of Rota--Baxter algebras and their relations to integral calculus and MZVs as well as for their general properties.

Let $A$ be a commutative Rota-Baxter algebra of weight $\theta$, that is a commutative algebra such that
$$\forall x,y\in A,\ R(x)R(y)=R(R(x)y+xR(y)+\theta xy).$$
The term $R(x)y+xR(y)+\theta xy=:x\ast_Ry$ defines a new commutative (and associative) product $\ast_R$ on $A$ called the \it double Rota-Baxter \rm product.
 We define the \it double quasi-shuffle bialgebra \rm over a Rota--Baxter algebra $A$, $\QSH^R(A)$, as the bialgebra which identifies with $T(A):=\bigoplus\limits_{n\in\N}A^{\otimes n}$ as a vector space, equipped with the deconcatenation coproduct, and equipped with the following recursively defined product $\shuffle_R$:
$$x_1...x_k\shuffle_R y_1...y_l:=x_1(x_2...x_k\shuffle_R y_1...y_l)+y_1(x_1...x_k\shuffle_R y_2...y_l)+$$
$$(x_1\ast_R y_1)(x_2...x_k\shuffle_R y_2...y_l).$$
The fact that $\QSH^R(A)$ is indeed a bialgebra follows from the general definition of the quasi-shuffle bialgebra over a commutative algebra $A$, see \cite{Hoffman,fpt}. 
\end{example}

\begin{example} This fourth example (a particular case of the previous one) and the following one are the first concrete examples of the kind of Hopf algebraic structure showing up specifically in mould calculus. The
definitions we introduce are inspired by the notion of \it symmetrul mould \rm \cite[p. 418]{ARI2003} of which they aim at capturing the underlying combinatorial structure. 

Let $\R [X]$ be equipped with the Riemann integral $R:=\int_0^X$ viewed as a Rota--Baxter operator of weight zero. With the notation $a_i:=X^{i-1},\ i\in\N^\ast$ we get: $R(a_i):=\frac{a_{i+1}}{i}$ and
$$a_i\ast_Ra_j=\frac{i+j}{ij}a_{i+j}.$$
This associative and commutative product gives rise to the following definition:
\begin{defn}
The bialgebra of quasi-symmetrul functions $\mathbf{QSul}$ is the quasi-shuffle bialgebra over the linear span of the integers $\N^\ast$ equipped with the product
$$[i]\ast [j]:=\frac{i+j}{ij}[i+j].$$
\end{defn}
\begin{prop}
The bialgebras $\QSym$, $\SH$ and  $\mathbf{QSul}$ are isomorphic, the isomorphism $\phi$ from $\SH$ to  $\mathbf{QSul}$
is given by:
$$\phi([n_1\dots n_k]):=\sum\limits_{\mu_1+\dots \mu_i=k}\frac{(n_1+\dots+n_{\mu_1})\dots (n_{\mu_1+\dots+\mu_{i-1}+1}+\dots+n_{\mu_1+\dots+\mu_{i}})}{\mu_1!\dots\mu_i!n_1\dots n_k}\cdot$$
$$[n_1+\dots+n_{\mu_1},\dots,n_{\mu_1+\dots+\mu_{i-1}+1}+\dots+n_{\mu_1+\dots+\mu_{i}}].$$
\end{prop}
The Theorem is an application of Hoffman's structure theorems for quasi-shuffle bialgebras \cite{Hoffman}. It also follows from the combinatorial analysis of quasi-shuffle bialgebras understood as deformations of shuffle bialgebras in \cite{fpt}.

\end{example}

\begin{example} The previous example gives the pattern for the notion of symmetrulity (and gives an hint for its analytic meaning).
Let now $V$ be a vector space with a distinghished basis ${\mathcal B}:=(v_i)_{i\in I}$ and $M$ a subsemigroup of $\R^{\ast,+}$, the strictly positive real numbers.
Let $A$ be the linear span of $M\times \mathcal B$ whose elements $(m,v)$ are represented $m\choose v$ to stick to the ``bimould'' calculus notation \cite{ARI2003}. We set:
$${m_1\choose v_1}\ast {m_2\choose v_2}:=-\frac{1}{m_2}{m_1+m_2\choose v_1}-\frac{1}{m_1}{m_1+m_2\choose v_2}.$$

Using the notation $m_1\ ...\ m_n\choose v_1\ ...\ v_n$ for the tensor product of the ${m_i\choose v_i}$ in $T(A)$, equipped with the deconcatenation coproduct,
the following recursively defined product defines a bialgebra structure denoted $\QSUL (M,V)$ on $T(A)$: 
$${m_1\ ...\ m_n\choose v_1\ ...\ v_n }\shuffle_{ul} {p_1\ ...\ p_k\choose w_1\ ...\ w_l}:={m_1\choose v_1}\left({m_2\ ...\ m_n\choose v_2\ ...\ v_n }\shuffle_{ul} {p_1\ ...\ p_k\choose w_1\ ...\ w_l}\right)
$$
$$+{p_1\choose w_1}\left(
{m_1\ ...\ m_n\choose v_1\ ...\ v_n }\shuffle_{ul} {p_2\ ...\ p_k\choose w_2\ ...\ w_l}\right)$$
$$
-\left(\frac{1}{m_2}{m_1+m_2\choose v_1}+\frac{1}{m_1}{m_1+m_2\choose v_2}\right)
{m_2\ ...\ m_n\choose v_2\ ...\ v_n }\shuffle_{ul} {p_2\ ...\ p_k\choose w_2\ ...\ w_l}.
$$

The property follows from the associativity and commutativity of the $\ast$ product, whose proof is left to the reader.
\end{example}

\begin{defc} For $B$ an arbitrary commutative algebra, writing $G_{\QSUL (M,V)}(B)$ for the group of $B$-valued characters of  $\QSUL (M,V)$, the functor $G_{\QSUL (M,V)}$ is a prounipotent affine group scheme whose points are called
symmetrul moulds \cite{ARI2003}.
\end{defc}

Let us write $\SH (M,V)$ for the shuffle bialgebra over $A$, and $G_{\SH (M,V)}$ for the corresponding prounipotent affine group scheme, we have:
\begin{thm}\label{isoqsul}
The prounipotent affine group schemes $G_{\QSUL (M,V)}$ and $G_{\SH (M,V)}$ are isomorphic. The isomorphism is induced by the bialgebra isomorphism $\phi$ between $\SH (M,V)$ and $\QSUL (M,V)$ defined by:
$$\phi{m_1\dots m_n\choose v_1\dots v_n}:=
\sum\limits_{\mu_1+\dots \mu_i=n}\frac{(-1)^{n-i}}{\mu_1!\dots\mu_i!}\cdot $$
$$\left(\sum\limits_{i=1}^{\mu_1}\frac{1}{m_1\dots m_{i-1}m_{i+1}\dots m_{\mu_1}}{m_1+\dots +m_{\mu_1}\choose v_i}\right)\dots$$
$$\left(\sum\limits_{i=\mu_1+\dots \mu_{i-1}+1}^{n}\frac{1}{m_{\mu_1+\dots \mu_{i-1}+1}\dots m_{i-1}m_{i+1}\dots m_n}{m_{\mu_1+\dots \mu_{i-1}+1}+\dots +m_{n}\choose v_i}\right).$$
\end{thm}
 
The Theorem follows once again from Hoffman's structure theorem for quasi-shuffle algebras by identification of the coefficients of the exponential isomorphism in the particular case under consideration.

In Ecalle's terminology, symmetrul moulds and symmetral moulds on $M\times V$ are canonically in bijection. Notice that whereas the definition of symmetrul moulds as characters on the Hopf algebra $\QSUL (M,V)$ is essentially a group-theoretical interpretation of the definitions given in \cite{ARI2003}, the equivalence between the two notions of symmetrulity and symmetrality of Thm \ref{isoqsul} (and therefore also the precise formula for the isomorphism) is new at our best knowledge.

We do not insist further on the notion of symmetrulity that is relatively easy to handle group-theoretically as we just have seen, and will focus preferably of the one of symmetrility whose signification for MZVs seems deeper and for which a group-theoretical account is harder to obtain since it does not seem possible to interpret symmetril moulds as elements of a prounipotent affine group scheme, but only as elements of a properly defined prounipotent group scheme.

\section{Generic bialgebras}

Symmetril moulds, of which a formal definition will be given later on, behave very much as characters on $\QSym$ or $T(x,y)$. There are even some conversion rules to move from one notion to the other, that we will explain later. Unfortunately, 
this notion of symmetrility fails to be accounted for by using a naive theory of characters on a suitable Hopf algebra. 
The aim of this section is to explain what has to be changed in the classical theory of Hopf algebras to make sense of the notion. 

The constructions in this section are motivated by the two notions of twisted bialgebras (also called Hopf species) explored in \cite{PR,PS1,PS2,AM} and the one of constructions in the sense of Eilenberg and MacLane \cite{P3}. However, both the theory of constructions and vector species are too functorial to account for the very specific combinatorics of symmetrility, and we have to introduce for its proper understanding a different framework. In view of the similarities with the theory of constructions, we decided to keep however the terminology of ``generic structures'' used in \cite{P3}.

Let $X$ be a finite or countable alphabet, partitioned into subsets $X=\coprod\limits_{i\in I}X_i$. We say that the partition is trivial if the $X_i$ are singletons.
A word over $X$ (possibly empty) is said to be \it generic \rm if it contains at most one letter in each $X_i$. If the partition is trivial, this means that no letter can appear twice. Similarly, tensor products of words are generic if they contain overall at most one letter in each $X_i$. Two generic tensor products of words, $w,w'$, are said to be \it in generic position \rm if $w\otimes w'$ is again a generic tensor product. Two linear combinations of generic tensor products of words $\sum_w\lambda_ww$, 
$\sum_{w'}\lambda_{w'}w'$ are in generic position if all the pairs $(w,w')$ are.
The underlying word $u(t)$ of a tensor product $t$ of words is the word obtained by concatenating its components: $u(x_1\otimes x_2x_4\otimes x_3)=u(x_1x_2\otimes x_4x_3)=x_1x_2x_4x_3$, so that $u(t)$ is generic if and only if $t$ is generic.

\begin{defn}
The category ${\mathbf {Gen}}_X$ of generic expressions over $X$ is the smallest linear (i.e. such that $Hom$-sets are $k$-vector spaces) subcategory of the category of vector spaces containing the null vector space and (examples refer to the case where $X=\{x_i\}_{i\in\N^\ast}$, with the trivial partition)
\begin{itemize}
\item containing the one-dimensional vector spaces $V_t$ generated by generic tensor products of words $t$,
\item closed by direct sums (although this won't be the case in the examples we will consider, multiple copies of the $V_t$ can be allowed, the following rules are applied to each of these copies)
\end{itemize}
and such that furthermore $Hom$ sets contain:
\begin{itemize}
\item for $t,t'$ two generic tensors with $u(t)=u(t')$, the map from $V_t$ to $V_{t'}$ induced by $f(t):=t'$,
\item the maps induced by substitutions of the letters inside the blocks $X_i$,
\item the maps obtained by erasing letters in the tensor products (e.g. the map induced by $f(x_1\otimes  x_2x_4\otimes x_3):=x_1\otimes   x_4$).
\end{itemize}
\end{defn}

Most importantly for our purposes, ${\mathbf {Gen}}_X$ is equipped with a symmetric monoidal category structure by the generic tensor product $\hat\otimes$ defined on the $V_t$ by 
$V_t\hat\otimes V_{t'}:= V_{t\otimes t'}$ if $t\otimes t'$ is generic, and $:=0$ else. The generic tensor product is extended to direct sums by the rule $(A\oplus B) \hat\otimes (C\oplus D)=A\hat\otimes C\oplus 
A\hat\otimes  D\oplus
 B \hat\otimes C\oplus
B \hat\otimes D$.
Notice, for further use, the canonical embedding $A\hat\otimes B\hookrightarrow A\otimes B$.

The reader familiar with homological algebra will have recognized the main ingredients of the theory of constructions \cite{P3}.
\it Generic algebras\rm , coalgebras, Hopf algebras, Lie algebras, and so on, are, by definition, algebras, coalgebras, Hopf algebras, Lie algebras, and so on, in a given ${\mathbf {Gen}}_X$.
For example, a generic algebra $A$ without unit is an object of ${\mathbf {Gen}}_X$ equipped with an associative product map $\mu$ from $A\hat\otimes A$ to $A$. Notice that $\mu$ can be viewed alternatively as a partially defined product map on $A$ (it is defined only on elements in $A\otimes A$ in generic position and linear combinations thereof). 

We will study from now on only \it standard \rm generic bialgebras $H$. By which it will be meant that $H$ is a generic bialgebra with product $\pi$ and coproduct $\Delta$ such that
\begin{itemize}
\item $H=\bigoplus\limits_{n\in\N}H_n$, where $H_0=V_\emptyset$ is identified with the ground field $k$ and $\emptyset$ behaves as a unit/counit for the product and the coproduct,
\item the coproduct is graded, 
\item the product satisfies the filtering condition: $\forall k,l>0,\ \pi(H_k\otimes H_l)\subset \bigoplus_{0<n\leq k+l}H_n$.
\end{itemize}
These bialgebras behave as the analogous usual bialgebras (the same arguments and proofs apply, we refer e.g. to \cite{Cartier} for the classical case). In particular
such a bialgebra is equipped with a convolution product of linear endomorphims: for arbitrary $f,g\in Hom_{{\mathbf{Gen}}_X}(H,H)$, $f\ast g:=\pi\circ (f\hat\otimes g)\circ \Delta$. The projection $u$ from $H$ to $H_0$ orthogonally to the $H_i,\ i\geq 1$ is a unit for $\ast$. 
Convolution of linear forms on $H$ is defined similarly.

The existence of an antipodal map, that is a convolution inverse $S$ 
to the identity map $I$ follows from the identity
\begin{equation}\label{antipode}
S=(u+(I-u))^{\ast \ -1}=u+\sum\limits_{n>0}(-1)^n(I-u)^{\ast n},
\end{equation}
where the rightmost sum restricts to a finite sum when $S$ is acting on a graded component $H_n$ since the coproduct is graded. In particular, a standard generic bialgebra $H$ is automatically a generic Hopf algebra.

Since $A\hat\otimes B\subset A\otimes B$, one can define morphisms from an algebra, bialgebra, Hopf algebra... in $Gen_X$ to a classical algebra, bialgebra, Hopf algebra... We will call such morphisms \it regularizing morphisms\rm .
For example, a regularizing morphism between a standard generic bialgebra $H$ equipped with the product $\mu$ and the coproduct $\Delta$ and a graded Hopf algebra $H'$ equipped with the product $\mu'$ and the coproduct $\Delta'$ is a morphism of graded vector spaces $f$ that maps the unit $\emptyset$ of $H$ to the unit $1\in H_0'$ of $H'$ and such that, for any $h,h'$ in generic position in $H$,
$$f(\mu(h\hat\otimes h'))=\mu'(f(h)\otimes f(h')), (f\otimes f)\circ \Delta (h)=\Delta'(f(h)).$$

\begin{example}
A first example of a standard generic bialgebra will look familiar to readers familiar with the theory of free Lie algebras and Reutenauer's monograph \cite{Reutenauer}. Let $X=[n]$ be equipped with the trivial partition. Then, let $T^g_k(X)$ be the linear span of generic words of length $k$: $T^g(X)=\bigoplus\limits_{n\in\N}T^g_n(X)$ is usually called the multilinear part of the tensor algebra over $X$ in the literature. Concatenation of words defines a map from $T^g_k(X)\hat\otimes T^g_l(X)$ to $T^g_{k+l}(X)$ and a generic algebra structure on $T^g(X)=\bigoplus\limits_{n\in\N}T^g_n(X)$. Similarly, the usual unshuffling of words $\Delta$ (the coproduct dual to the one introduced in example 2) defines, when restricted to generic words, a generic coalgebra structure, and, together with the concatenation product, a standard generic bialgebra structure on $T^g(X)$. The generic Lie algebra of primitive elements of $T^g(X)$, which is defined as usual:
$Prim(T^g(X)):=\{w\in T^g(X),\Delta (w)=w\otimes 1+1\otimes w\}$, is simply the multilinear part of the usual free Lie algebra over $X$.

Dually, the shuffle product and the deconcatenation product (as in Example \ref{ex2}) define a (dual) standard generic bialgebra structure on $T^g(X)$, that will be named the generic shuffle bialgebra over $X$ and denoted $\SH^g(X)$. We write simply $\SH^g$ for $\SH^g(\N^\ast)$
\end{example}

This example is particularly easy to understand: the embedding of $T^g(X)$ into the usual tensor algebra over $X$, $T(X)$, is a regularizing morphism, and all our assertions are direct consequences of the behaviour of $T(X)$ as exponed e.g. in \cite{Reutenauer}. 

\begin{defn}
 Let $H$ be a standard generic bialgebra. For $B$ an arbitrary commutative algebra, a $B$-valued character on $H$ is, by definition, a unital multiplicative map from $H$ to $B$, that is a map $\phi$ such that:
\begin{itemize}
\item $\phi(\emptyset)=1$,
\item For any $h_1,h_2$ in generic position, writing $h_1\cdot h_2:=\pi (h_1\hat\otimes h_2)$

\begin{equation}\label{condit}\phi (h_1\cdot h_2)=\phi(h_1)\phi(h_2).\end{equation}
\end{itemize}
\end{defn}

 \begin{prop}Let $H$ be a standard generic bialgebra.
The set $G_H(B)$ of $B$-valued characters is equipped with a group structure by the convolution product $\ast$. The corresponding functor $G_H$ from commutative algebras over the reference ground field $k$ to groups is called, by analogy with the classical case, a generic group scheme.
\end{prop}
 Indeed, we have, for any $\phi,\phi'\in G_H(B)$, and any $h,h'\in H^+:=\bigoplus\limits_{n>0}H_n$ in generic position:
 $$\phi\ast\phi'(\emptyset)=\phi(\emptyset)\phi'(\emptyset)=1,$$
 $$\phi\ast\phi'(h\cdot h')=\phi(h^{(1)}\cdot h^{'(1)})\phi'(h^{(2)}\cdot h^{'(2)})$$
 $$=\phi(h^{(1)})\phi'(h^{(2)})\phi(h^{'(1)})\phi'(h^{'(2)})$$
 $$=\phi\ast\phi'(h)\cdot \phi\ast\phi'(h'),$$
  where we used a Sweedler-type notation $\Delta(h)=h^{(1)}\hat\otimes h^{(2)}$.
  
  Similarly, $\phi\circ S$ is the convolution inverse of $\phi$ since:
  $((\phi\circ S)\ast\phi)(\emptyset)=\phi(\emptyset)^2=1$ and, for $h$ as above,
  $$((\phi\circ S)\ast\phi )(h)=\phi(S(h^{(1)}))\phi(h^{(2)})=\phi(S(h^{(1)})\cdot h^{(2)})=\phi\circ u(h)=0.$$
  Notice that, contrary to the classical case, the identity $\phi(S(h^{(1)}))\phi(h^{(2)})=\phi(S(h^{(1)})\cdot h^{(2)})$ is not straightforward since identity \ref{condit} holds only under the asumption that $h_1,h_2$ are in generic position. Here, we can apply the identity  because $S$, in view of Eq. \ref{antipode}, can be written on each graded component as a sum of convolution powers $I^{\ast k}$ of the identity map. It is then enough to check that, given $h\in H^+$, $I^{\ast k}(h^{(1)})\otimes h^{(2)}$ can be written as a linear combination of tensor products $w\otimes w'$, where $w,w'$ are in generic position, which follows from the definition of the convolution product $\ast$ and the coassociativity of $\Delta$.

\section{Symmetril moulds and generic group schemes}
We come now to the main examples of generic structures in view of the scope of the present article --symmetrility properties. This section aims at abstracting the key combinatorial features of symmetrility in order to study them and link them with classical combinatorial objects, such as quasi-symmetric functions. The next section will move forward by sticking closer to Ecalle's study of MZVs, linking symmetrility phenomena to the resummation of MZVs.

\begin{defn}
Let $X=\N^\ast$, equipped with the trivial partition. We define the generic divided quasi-shuffle bialgebra over $\N^\ast$, $\QSH^g_d$, as the generic bialgebra which identifies with $T^g(\N^\ast)$ as a vector space, equipped with the deconcatenation coproduct, and equipped with the following recursively defined product $\dshuffle$ (elements of $T^g(\N^\ast)$ are written using a bracketed word notation):
$$[n_1...n_k]\dshuffle [m_1...m_l]:=[n_1(n_2...n_k\dshuffle m_1...m_l)]+[m_1(n_1...n_k\dshuffle m_2...m_l)]+$$
$$\frac{1}{n_1-m_1}\{[n_1(n_2...n_k\dshuffle m_2...m_l)]-[m_1(n_2...n_k\dshuffle m_2...m_l)]\}.$$
The elements of the groups $G_{\QSH^g_d}(B)$ are called symmetril moulds (over $\N^\ast$).
\end{defn}

Proving that $\QSH^g_d$ is indeed a Hopf algebra in ${\mathbf {Gen}}_X$ is not entirely straightforward and is better stated at a more general level, by mimicking for generic structures the theory of quasi-shuffle algebras.

\begin{defp}
Let $X$ be a partitioned alphabet and assume that $\ast$ equips $k<X>$, the linear span of $X$, with the structure of a generic commutative algebra. Then, the generic quasi-shuffle bialgebra denoted $\QSH_\ast^g(X)$ over $(k<X>,\ast)$ is, by definition, the generic bialgebra whose underlying generic coalgebra is $T^g(X)$ equipped with the deconcatenation coproduct $\Delta$, and whose commutative  product is defined inductively (for words satisfying the genericity conditions) by:
$$[n_1...n_k]\dshuffle [m_1...m_l]:=[n_1(n_2...n_k\dshuffle m_1...m_l)]+[m_1(n_1...n_k\dshuffle m_2...m_l)]+$$
$$[(n_1\ast m_1)(n_2...n_k\dshuffle m_2...m_l) ].$$
\end{defp}

The associativity of the product follows by induction on the total length $k+l+q$ from the identity of the expansion:
$$[(n_1...n_k\dshuffle m_1...m_l)\dshuffle p_1..p_q]=[
n_1(n_2...n_k\dshuffle m_1...m_l\dshuffle p_1..p_q)]$$
$$+[p_1((n_1(n_2...n_k\dshuffle m_1...m_l))\dshuffle p_2..p_q)]+[(n_1\ast p_1)(n_2...n_k\dshuffle m_1...m_l\dshuffle p_1..p_q)]$$
$$+
[m_1(n_1...n_k\dshuffle m_2...m_l\dshuffle p_1..p_q)]+[p_1((m_1(n_1...n_k\dshuffle m_2...m_l))\dshuffle p_2..p_q)]$$
$$+[
(m_1\ast p_1)(n_1...n_k\dshuffle m_2...m_l\dshuffle p_2...p_q)]+[(n_1\ast m_1)(n_2...n_k\dshuffle m_2...m_l\dshuffle p_1..p_q)] +$$
$$+[p_1((n_1\ast m_1)(n_2...n_k\dshuffle m_2...m_l))\dshuffle p_2..p_q)]+[(n_1\ast m_1\ast p_1)(n_1...n_k\dshuffle m_2...m_l\dshuffle p_1...p_q)]$$
$$=[n_1(n_2...n_k\dshuffle m_1...m_l\dshuffle p_1..p_q)]+[m_1(n_1...n_k\dshuffle m_2...m_l\dshuffle p_1..p_q)]$$
$$+[p_1(n_1...n_k\dshuffle m_1...m_l\dshuffle p_2...p_q)]+[(n_1\ast m_1)(n_2...n_k\dshuffle m_2...m_l\dshuffle p_1..p_q)] $$
$$
+[(n_1\ast p_1)(n_2...n_k\dshuffle m_1...m_l\dshuffle p_2...p_q)]+
[(m_1\ast p_1)(n_1...n_k\dshuffle m_2...m_l\dshuffle p_2...p_q)]
$$
$$
+[(n_1\ast m_1\ast p_1)(n_1...n_k\dshuffle m_2...m_l\dshuffle p_2...p_q)]$$
with the same symmetric expansion in the $n_i,m_i,p_i$ for $[n_1...n_k\dshuffle (m_1...m_l\dshuffle p_1...p_q)]$.

The compatibility of the deconcatenation coproduct with the product is obtained similarly and follows the same pattern than the proof that usual quasi-shuffle algebras over commutative algebras are indeed equipped with a Hopf algebra structure by the deconcatenation coproduct \cite{Hoffman,fpt}, and is omitted.

We can now conclude that $\QSH_d^g$ is indeed a generic  bialgebra from the Lemma:

\begin{lem}
The product $\ast$,
$$[n]\ast [m]:=\frac{1}{n-m}([n]-[m])$$
equips $k<\mathbf N^\ast>$ with the structure of a generic commutative algebra.
\end{lem}

Indeed, for distinct $m,n,p$,
$$[(n\ast m)\ast p]=\frac{1}{n-m}([n]-[m])\ast [p]=\frac{1}{(n-m)(n-p)}[n]$$
$$+\frac{1}{(m-n)(m-p)}[m]+(\frac{1}{(n-m)(p-n)}+\frac{1}{(m-n)(m-p)})[p]$$
$$=\frac{1}{(n-m)(n-p)}[n]+\frac{1}{(m-n)(m-p)}[m]+\frac{1}{(p-n)(p-m)}[p]$$
which is equal to the same symmetric expression for $[n\ast (m\ast p)]$.

For later use, we also calculate iterated products in $k<\mathbf N^\ast>$.
\begin{lem}For distinct $n_1,\dots,n_k\in N^\ast$ we have
$$[n_1]\ast \dots\ast [n_k]:=\sum\limits_{i=1}^k\frac{[n_i]}{\prod\limits_{j\not= i}(n_i-n_j)}$$
\end{lem}

Let us assume that the Lemma holds for $k\leq p$ and prove it by induction. Since the product $\ast$ is commutative, it is enough to show that the coefficient of $ [n_{p+1}]$ in $[n_1]\ast \dots\ast [n_{p+1}]$ is given by 
$\frac{1}{\prod\limits_{j\leq p}(n_{p+1}-n_j)}$.
Equivalently, we have to show that $\alpha =1$, where
$$\alpha=\sum\limits_{i=1}^p\frac{\prod\limits_{j\leq p}(n_{p+1}-n_j)}{(n_{p+1}-n_i)\prod\limits_{j\not= i,j\leq p}(n_i-n_j)}=\sum\limits_{i=1}^p\prod\limits_{j\not= i, j\leq p}\frac{(n_{p+1}-n_j)}{(n_i-n_j)}.$$

Notice that the induction hypothesis amounts to assuming that the following two equivalent identities hold for arbitrary distinct integers $m_1,...,m_p$
$$\sum\limits_{i=1}^{p-1}\prod\limits_{j\not= i,j\leq p-1}\frac{(m_p-m_j)}{(m_i-m_j)}=1,\ \sum\limits_{i=1}^{p}\prod\limits_{j\not= i,j\leq p}\frac{1}{(m_i-m_j)}=0.$$

We get:
$$\alpha=\sum\limits_{i=1}^{p-1} \left(\prod\limits_{j\not= i,j\leq p}\frac{(n_{p+1}-n_j)}{(n_i-n_j)}\right)+\prod\limits_{j\leq p-1}\frac{(n_{p+1}-n_j)}{(n_p-n_j)}$$

$$=\sum\limits_{i=1}^{p-1} \left(\frac{\prod\limits_{j\not= i,j\leq p-1}(n_{p+1}-n_j)}{\prod\limits_{j\not= i,j\leq p}(n_i-n_j)}\right)\left( (n_{p+1}-n_i)+(n_i-n_p)\right)
+\prod\limits_{j\leq p-1}\frac{(n_{p+1}-n_j)}{(n_p-n_j)}$$
$$=\sum\limits_{i=1}^{p} \left(\frac{1}{\prod\limits_{j\not= i,j\leq p}(n_i-n_j)}\right)
\cdot \prod\limits_{j\leq p-1}(n_{p+1}-n_j)+
\sum\limits_{i=1}^{p-1} \left(\prod\limits_{j\not= i,j\leq p-1}\frac{(n_{p+1}-n_j)}{(n_i-n_j)}\right)=0+1=1,$$
where the last identity follows from the induction hypothesis.

\begin{thm}\label{regSH}
The following map $\psi$ defines a linear embedding of $\QSH_d^g$ into $\SH$ and is a regularizing bialgebra map. 
$$\psi([n_1\dots n_k]):=\sum\limits_{\mu_1+\dots+\mu_i=k}\frac{(-1)^{k-i}}{\mu_1\dots\mu_i}\left(\sum\limits_{j=1}^{\mu_1}\frac{[n_j]}{\prod\limits_{l\not= j,l\leq \mu_1}(n_j-n_l)}\right) \dots$$
$$\dots \left(\sum\limits_{j=\mu_1+\dots+\mu_{i-1}+1}^{k}\frac{[n_j]}{\prod\limits_{l\not= j,\mu_1+\dots+\mu_{i-1}+1\leq l\leq k}(n_j-n_l)}\right) $$
In particular, the product and coproduct maps on $\QSH_d^g$ are mapped to the product and coproduct on $\SH$.
\end{thm}

The Theorem can be rephrased internally to the category ${\mathbf {Gen}}_{\N^\ast}$ -this is because the image of $\psi$ identifies with the subspace $T^g(\N^\ast)$ of $\SH$ (the latter identifying with $T(\N^\ast)$ as a graded vector space).

\begin{cor}\label{isogen}
The standard generic bialgebras $\QSH_d^g$ and $\SH^g$ are isomorphic under $\psi$. 
\end{cor}

The theorem is an extension to the generic case of the Hoffman isomorphism between shuffle and quasi-shuffle bialgebras. Following \cite{Hoffman,fpt}, the proof of the isomorphism relies only on the combinatorics of partitions and on a suitable lift to formal power series of natural coalgebra endomorphisms of shuffle bialgebras (we refer to \cite{fpt} for details). Let us show here that these arguments still hold in the generic framework.

Let $P(X)=\sum\limits_{i=1}^\infty p_iX^i$ be a formal power series $X \mathbf Q[[X]]$.
This power series induces a generic coalgebra endomorphism $\phi_{P}$
of $T^g(\N^\ast)$ equipped with the deconcatenation product:
 on an arbitrary generic tensor $[n_1\dots n_k]\in T^g(\N^\ast)$ the action is given by
 \begin{equation}\label{flaphi}
\phi_{P}([n_1\dots n_k])=\sum\limits_{j=1}^k\sum\limits_{i_1+...+i_j=k}p_{i_1}...p_{i_j}([n_1]\ast ...\ast [n_{i_1}])\otimes ...\otimes ([n_{i_1+...+i_{j-1}+1}]\ast ...\ast [n_{k}]),
\end{equation}
where we recall that $[n]\ast[m]:=\frac{[n]-[m]}{n-m}$. 
When $p_1\not=0$, $\phi_P$ is bijective (by a triangularity argument), and a coalgebra automorphism of $T^g(\N^\ast)$. 

Let us show now that, for arbitrary $P(X),\ Q(X)\in X\mathbf Q[[X]]$, 
\begin{equation}
\phi_{P}\circ\phi_{Q}=\phi_{P\circ Q},
\end{equation}
where $(P\circ Q)(X):=P(Q(X))$.
We have indeed, for an arbitrary sequence of distinct integers $n_1,\dots,n_k$:
$$\phi_{P}\circ\phi_{Q}(n_1...n_k)=$$
$$=\phi_{P}(\sum\limits_{j=1}^k\sum\limits_{i_1+\dots +i_j=k}q_{i_1}\dots q_{i_k}(n_1\ast\dots\ast n_{i_1})\otimes\dots\otimes (n_{i_1+\dots+i_{j-1}+1}\ast \dots\ast n_{k}))$$
$$=\sum\limits_{j=1}^k\sum\limits_{l=1}^j\sum\limits_{h_1+\dots +h_l=j}\sum\limits_{i_1+\dots +i_j=k}p_{h_1}\dots p_{h_l}q_{i_1}\dots q_{i_k}(n_1\ast\dots\ast n_{i_1+\dots+i_{h_1}})
\otimes$$
$$\dots\otimes (n_{i_1}+\dots+i_{h_1+\dots+h_{l-1}+1}\ast \dots\ast n_{k})$$
$$=\phi_{P(Q)}(n_1...n_k).$$
The proof of the theorem follows: $\psi=\phi_{log}$ has for inverse $\rho=\phi_{exp}$ which maps isomorphically $\SH^g$ to $\QSH_d^g$ (Hoffman's combinatorial argument in the classical case in \cite{Hoffman} applies \it mutatis mutandis \rm when restricted to generic tensors).

%%%%%%%%%%%%%%%%%%%%%%%%%%%%%%%%%%%%%%%%%

\section{Resummation of MZVs}
In order to resum MZVs into formal power series equipped with interesting group-theoretical operations and structures, let us introduce first a formal analog of the standard generic bialgebra $\QSH_d^g$ studied previsouly. Here, "formal" means that numbers and sequences of numbers are replaced by formal power series and words over an alphabet. Proofs of the properties and structure theorems are similar to the ones for $\QSH_d^g$ and are omitted. Our definitions and constructions are motivated by \cite{ARI2003}.

\begin{defn}
Let $V=\{v_i\}_{i\in\N^\ast}$, equipped with the trivial partition. We define the generic divided quasi-shuffle bialgebra over $V$, $\QSH_d^g(V)$, as the generic bialgebra defined over $k_V:=k((V))$, the field of fractions of the ring of formal power series over $V$, which identifies with $T^g(V)$ as a vector space, equipped with the deconcatenation coproduct, and equipped with the following recursively defined product $\dshuffle$ (elements of $T^g(V)$ are written using a bracketed word notation):
$$[v_{i_1}...v_{i_k}]\dshuffle [v_{i_{k+1}}...v_{i_{k+l}}]:=[v_{i_1}(v_{i_2}...v_{i_k}\dshuffle v_{i_{k+1}}...v_{i_{k+l}})]+[v_{i_{k+1}}(v_{i_1}...v_{i_k}\dshuffle v_{i_{k+2}}...v_{i_{k+l}})]+$$
$$\frac{1}{v_{i_1}-v_{i_{k+1}}}\{[v_{i_1}(v_{i_2}...v_{i_k}\dshuffle v_{i_{k+2}}...v_{i_{k+l}})]-[v_{i_{k+1}}(v_{i_2}...v_{i_k}\dshuffle v_{i_{k+2}}...v_{i_{k+l}})]\},$$
where $[v_{i_1}...v_{i_k}]$ and $[v_{i_{k+1}}...v_{i_{k+l}}]$ are in generic position (so that $\frac{1}{v_{i_1}-v_{i_{k+1}}}$ is well-defined).

The elements of the groups $G_{\QSH^g_d(V)}(B)$, where $B$ runs over algebras over $k_V$, associated to the generic group scheme $G_{\QSH^g_d(V)}$ over $k_V$ are called symmetril moulds (over $V$).
\end{defn}

Let us denote $\SH_V^g$ the generic shuffle bialgebra over $V$ with $k_V$ as a field of coefficients. Corollary \ref{isogen} generalizes to $\QSH_d^g(V)$ and $\SH_V^g$: the two generic bialgebras are isomorphic under $\psi_V$:
$$\psi_V([v_1\dots v_k]):=\sum\limits_{\mu_1+\dots+\mu_i=k}\frac{(-1)^{k-i}}{\mu_1\dots\mu_i}\left(\sum\limits_{j=1}^{\mu_1}\frac{[v_j]}{\prod\limits_{l\not= j,l\leq \mu_1}(v_j-v_l)}\right) \dots$$
$$\dots \left(\sum\limits_{j=\mu_1+\dots+\mu_{i-1}+1}^{k}\frac{[v_j]}{\prod\limits_{l\not= j,\mu_1+\dots+\mu_{i-1}+1\leq l\leq k}(v_j-v_l)}\right) .$$

Let us denote now $\QSym_V$ the completion (w.r. to the grading) of the bialgebra of quasi-symmetric functions over the base field $k_V$. Since properly regularized MZVs at positive values are characters on $\QSym$, generating series for MZVs such as $\sum\limits_{n_1,...,n_k\geq 1}v_1^{n_1-1}\dots v_k^{n_k-1}\zeta(n_1,\dots,n_k)$ and the study of their algebraic structure can be lifted to $\QSym_V$. Let us show how this idea translates group-theoretically.

\begin{thm}
The following morphism $\gamma$ is a regularizing bialgebra map from $\QSH^g_d(V)$ to $\QSym_V$:
$$\gamma([v_{i_1}\dots v_{i_k}]):=\sum\limits_{n_1,\dots,n_k\geq 1}v_{i_1}^{n_1-1}\dots v_{i_k}^{n_k-1}\cdot [n_1\dots n_k].$$
\end{thm}

Notice first that $\gamma$ is, by its very definition, multiplicative for the concatenation product:
\begin{equation}\label{concat}
\gamma([v_{i_1}\dots v_{i_k}])=\gamma([v_{i_1}])\cdot \gamma([v_{i_2}\dots v_{i_k}])=\gamma([v_{i_1}])\cdot \gamma([v_{i_2}])\dots\gamma([v_{i_k}]),
\end{equation}
from which it follows that $\gamma$ is a coalgebra map (recall that the later is induced on $\QSH^g_d(V)$ and $\QSym_V$ by deconcatenation).

Let us prove that, for any $v_{i_1}...v_{i_k}$, $v_{i_{k+1}}...v_{i_{k+l}}$ in generic position, we have
$$\gamma([v_{i_1}\dots v_{i_k}]\dshuffle [v_{i_{k+1}}\dots v_{i_{k+l}}] )=
\gamma([v_{i_1}\dots v_{i_k}])\qshuffle \gamma([v_{i_{k+1}}\dots v_{i_{k+l}}] )$$
by induction on $k+l$. So, let $v_{i_0}$ an element of $V$ distinct from $v_{i_1},...,v_{i_{k+l}}$. We get:
$$\gamma([v_{i_0}\dots v_{i_k}]\dshuffle [v_{i_{k+1}}\dots v_{i_{k+l}}] )=
\gamma([v_{i_0}(v_{i_1}...v_{i_k}\dshuffle v_{i_{k+1}}...v_{i_{k+l}})])$$
$$+\gamma([v_{i_{k+1}}(v_{i_0}...v_{i_k}\dshuffle v_{i_{k+2}}...v_{i_{k+l}})])$$
$$+\gamma(\frac{1}{v_{i_0}-v_{i_{k+1}}}\{[v_{i_0}(v_{i_1}...v_{i_k}\dshuffle v_{i_{k+2}}...v_{i_{k+l}})]-[v_{i_{k+1}}(v_{i_1}...v_{i_k}\dshuffle v_{i_{k+2}}...v_{i_{k+l}})]\}).$$
From Eqn \ref{concat} and the induction hypothesis, we get:
$$\gamma([v_{i_0}(v_{i_1}...v_{i_k}\dshuffle v_{i_{k+1}}...v_{i_{k+l}})])=
\gamma([v_{i_0}]) \gamma([v_{i_1}...v_{i_k}\dshuffle v_{i_{k+1}}...v_{i_{k+l}})])$$
$$=\gamma([v_{i_0}])\gamma([v_{i_1}\dots v_{i_k}])\qshuffle \gamma([v_{i_{k+1}}\dots v_{i_{k+l}}] )$$
and similarly
$$\gamma([v_{i_{k+1}}(v_{i_0}...v_{i_k}\dshuffle v_{i_{k+2}}...v_{i_{k+l}})])=
\gamma([v_{i_{k+1}}])(\gamma([v_{i_0}\dots v_{i_k}])\qshuffle \gamma([v_{i_{k+2}}\dots v_{i_{k+l}}] )).$$
At last,
$$\gamma(\frac{1}{v_{i_0}-v_{i_{k+1}}}[v_{i_0}-v_{i_{k+1}}][v_{i_1}...v_{i_k}\dshuffle v_{i_{k+2}}...v_{i_{k+l}}])=$$
$$
\frac{1}{v_{i_0}-v_{i_{k+1}}}
\gamma([v_{i_0}-v_{i_{k+1}}])\gamma([v_{i_1}...v_{i_k}\dshuffle v_{i_{k+2}}...v_{i_{k+l}}])$$
and, in view of the recursive definition of $\qshuffle$, to conclude the proof it remains to show that
$$\frac{1}{v_{i_0}-v_{i_{k+1}}}
\gamma([v_{i_0}-v_{i_{k+1}}])=\gamma([v_{i_0}])\odot\gamma([v_{i_{k+1}}])$$
where, to avoid confusion with other already introduced symbols, $\odot$ denotes the product of bracketed integers induced by the addition: $[n]\odot[m]=[n+m]$.

Indeed, we have:
$$\gamma([v_{i_0}-v_{i_{k+1}}])=\sum_{n\geq 1}(v_{i_0}^{n-1}-v_{i_{k+1}}^{n-1})[n]=\sum_{n\geq 2}(v_{i_0}^{n-1}-v_{i_{k+1}}^{n-1})[n],$$
and
$$(v_{i_0}-v_{i_{k+1}})\gamma([v_{i_0}])\odot\gamma([v_{i_{k+1}}])=
(v_{i_0}-v_{i_{k+1}})\sum_{n,m\geq 1}v_{i_0}^{n-1}v_{i_{k+1}}^{m-1}[n+m]$$
$$=\sum_{p\geq 2}(v_{i_0}-v_{i_{k+1}})(\sum_{n,m\geq 0,n+m=p-2}v_{i_0}^{n}v_{i_{k+1}}^{m})[p]=\sum_{p\geq 2}(v_{i_0}^{p-1}-v_{i_{k+1}}^{p-1})[p].$$

\begin{cor}Let $V$ be an infinite alphabet.
The regularizing morphism $\gamma$ induces, for any commutative algebra $B$ over a base field $k$ a group map from $G_{\QSym}(B)$ to $G_{\QSH_d^g(V)}(B\otimes_kk_V)$.
\end{cor}

In particular, regularized $\zeta$ functions, viewed as a real-valued characters on $\QSym$ give rise to symmetril $\R((V))$-valued moulds over $V$. More generally, symmetrel moulds give rise to symmetril moulds by resummation \cite{ARI2003,cresson}-the very reason for the introduction of the latter.

Let us mention that, when dealing with modular MZVs,
$$\zeta {\epsilon_1\dots \epsilon_r\choose s_1\dots s_r}:=\sum\limits_{n_1>\dots>n_r}\frac{e^{2\pi in_1\epsilon_1}\dots e^{2\pi in_r\epsilon_r}}{n_1^{s_1}\dots n_r^{s_r}},$$
where $\epsilon_i\in{\mathbf Q}/\mathbf Z$, a ``bimould'' version of the previous construction has to be used. We only sketch the constructions in that case, they could be developed in more detail following the previous ones in this section.

\begin{defn}
Let $W:={\mathbf Q}/{\mathbf Z}\times V$, with $V
=\{v_i\}_{i\in\N^\ast}$, equipped with the partition $W=\coprod W_i$, $W_i:={\mathbf Q}/{\mathbf Z}\times \{v_i\}$. We define the generic divided quasi-shuffle bialgebra over $W$, $\QSH^g_d(W)$, as the generic bialgebra defined over $k_V:=k((V))$, which identifies with $T^g(W)$ as a vector space, equipped with the deconcatenation coproduct, and equipped with the following recursively defined product $\dshuffle$ (elements of $W$ are represented as column vector):
$${\epsilon_1\dots \epsilon_r\choose v_{i_1}\dots v_{i_r}}
\dshuffle {\epsilon_{r+1}\dots \epsilon_{r+s}\choose v_{i_{r+1}}\dots v_{i_{r+s}}}
:={\epsilon_1\choose v_{i_1}}\left({\epsilon_2\dots \epsilon_r\choose v_{i_2}\dots v_{i_r}}
\dshuffle {\epsilon_{r+1}\dots \epsilon_{r+s}\choose v_{i_{r+1}}\dots v_{i_{r+s}}}\right)$$
$$+{\epsilon_{r+1}\choose v_{i_{r+1}}}\left({\epsilon_1\dots \epsilon_r\choose v_{i_1}\dots v_{i_r}}
\dshuffle {\epsilon_{r+2}\dots \epsilon_{r+s}\choose v_{i_{r+2}}\dots v_{i_{r+s}}}\right)$$
$$+\frac{1}{v_{i_1}-v_{i_{r+1}}}{\epsilon_1+\epsilon_{r+1}\choose v_{i_1}}\left({\epsilon_2\dots \epsilon_r\choose v_{i_2}\dots v_{i_r}}
\dshuffle {\epsilon_{r+2}\dots \epsilon_{r+s}\choose v_{i_{r+2}}\dots v_{i_{r+s}}}\right)$$
$$-\frac{1}{v_{i_1}-v_{i_{r+1}}}{\epsilon_1+\epsilon_{r+1}\choose v_{i_{r+1}}}\left({\epsilon_2\dots \epsilon_r\choose v_{i_2}\dots v_{i_r}}
\dshuffle {\epsilon_{r+2}\dots \epsilon_{r+s}\choose v_{i_{r+2}}\dots v_{i_{r+s}}}\right)$$
where ${\epsilon_1\dots \epsilon_r\choose v_{i_1}\dots v_{i_r}}$ and $
 {\epsilon_{r+1}\dots \epsilon_{r+s}\choose v_{i_{r+1}}\dots v_{i_{r+s}}}$ are in generic position (so that $\frac{1}{v_{i_1}-v_{i_{k+1}}}$ is well-defined).

The elements of the groups $G_{\QSH^g_d(W)}(B)$, where $B$ runs over algebras over $k_V$, associated to the generic group scheme $G_{\QSH^g_d(W)}$ over $k_W$ are called symmetril moulds (over $W$).
\end{defn}

Symmetril moulds over $W$ can be used to resum modular MZVs by the same process that allows the resummation of usual MZVs by symmetril moulds over $V$, see \cite{ARI2003}.

\section{A new resummation process}

In this last section, we introduce a new resummation process for MZVs, based on Thm \ref{regSH}. Contrary to Ecalle's resummation process, which maps a symmetrel mould (a character on the algebra of quasi-symmetric functions) to a symmetril mould, the new resummation is much more satisfactory in that it maps a symmetrel mould to a character on $\SH^g_V$, so that calculus on MZVs and other characters on $\QSym$ can be interpreted in terms of the usual rules of Lie calculus (recall that the set of primitive elements in the dual of $\Sh_V^g$ is simply the multilinear part of the free Lie algebra over the integers, a well-known object whose study is even easier than the one of the usual free Lie algebra).

\begin{thm}The inverse $\rho_V$ of the standard generic bialgebra isomorphism $\psi_V$ between $\QSH_d^g (V)$ and $\SH_V^g$ is given by
$$\rho_V([v_1\dots v_k]):=\sum\limits_{\mu_1+\dots+\mu_i=k}\frac{1}{\mu_1!\dots\mu_i!}\left(\sum\limits_{j=1}^{\mu_1}\frac{[v_j]}{\prod\limits_{l\not= j,l\leq \mu_1}(v_j-v_l)}\right) \dots$$
$$\dots \left(\sum\limits_{j=\mu_1+\dots+\mu_{i-1}+1}^{k}\frac{[v_j]}{\prod\limits_{l\not= j,\mu_1+\dots+\mu_{i-1}+1\leq l\leq k}(v_j-v_l)}\right). $$
\end{thm}

The Theorem follows by adapting to $T^g(V)$ the correspondence between formal power series in $X\bQ [[X]]$ and generic coalgebra endomorphisms of $T^g(\bN^\ast)$: with the same notation than the ones used for $T^g(\bN^\ast)$, each $P\in X\bQ [[X]]$ defines a generic coalgebra endomorphism $\phi_P$  of $T^g(V)$. We have $\rho_V=\phi_{exp}$ and $\psi_V=\phi_{log}$, and the two morphisms are mutually inverse.

\begin{cor}
The morphism $reg_V:=\gamma\circ\rho_V$ is a regularizing Hopf algebra morphism from $\SH_V^g$ to $\QSym_V$.
It induces, for any commutative algebra $B$ over the base field $k$ a group map from $G_\QSym(B)$ to $G_{\SH_V^g}(B\otimes_kk_V)$. 
\end{cor}

Naming generic symmetral moulds the characters on $\SH_V^g$, we get that 
this last map resums symmetrel moulds (such as regularized MZVs at the positive integers) into generic symmetral moulds. As announced, this approach should provide a new way to investigate group-theoretically the properties of MZVs. Together with the study of the various combinatorial structures introduced in the present article, this will the object of further studies.

We conclude by illustrating the resummation process on low dimensional examples that illustrate the behaviour of the map $reg_V$. We write $\zeta$ for a character on $\QSym$ (a symmetrel mould), having in mind the example of regularized multizetas.
The morphism $reg_V$ is given in low degrees by:
$$reg_V([v_1])=\gamma([v_1])=\sum\limits_{n\geq 1}v_1^{n-1}[n],$$
$$reg_V([v_1,v_2])=\gamma([v_1,v_2]+\frac{1}{2}\frac{[v_1]-[v_2]}{v_1-v_2})$$
$$=\sum\limits_{n,m\geq 1}v_1^{n-1}v_2^{m-1}[n,m]+
\frac{1}{2(v_1-v_2)}\sum\limits_{n\geq 1}(v_1^{n-1}-v_2^{n-1})[n]$$
$$=\sum\limits_{n,m\geq 1}v_1^{n-1}v_2^{m-1}[n,m]+
\frac{1}{2}\sum\limits_{n\geq 2, p+q=n-2}v_1^{p}v_2^{q}[n].$$

$$reg_V([v_1,v_2,v_3])=\gamma([v_1,v_2,v_3]+\frac{1}{2}\left(\frac{[v_1v_3]-[v_2v_3]}{v_1-v_2}+\frac{[v_1v_2]-[v_1v_3]}{v_2-v_3}\right)$$
$$+\frac{1}{6}\left(\frac{[v_1]}{(v_1-v_2)(v_1-v_3)}+\frac{[v_2]}{(v_2-v_1)(v_2-v_3)}+\frac{[v_3]}{(v_3-v_1)(v_3-v_2)}\right)
$$
$$=\sum\limits_{n,m,p\geq 1}v_1^{n-1}v_2^{m-1}v_3^{p-1}[n,m,p]+
\frac{1}{2}(\sum\limits_{n\geq 2, p+q=n-2,m\geq 1}v_1^{p}v_2^{q}v_3^{m-1}[n,m]+$$
$$\sum\limits_{n\geq 1,m\geq 2, p+q=m-2}v_1^{n-1}v_2^{p}v_3^{q}[n,m]
)+\frac{1}{6}\sum\limits_{n\geq 3,\ p+q+r=n-3}v_1^{p}v_2^qv_3^r[n],$$
where we used the identity
$$\frac{v_1^{n-1}}{(v_1-v_2)(v_1-v_3)}+\frac{v_2^{n-1}}{(v_2-v_1)(v_2-v_3)}+\frac{v_3^{n-1}}{(v_3-v_1)(v_3-v_2)}$$
$$=\sum\limits_{p+q+r=n-3}v_1^{p}v_2^qv_3^r.$$

We get, for the $\zeta$ character:
$$\zeta\circ reg_V([v_1]\shuffle [v_2])=\zeta (\sum\limits_{n,m\geq 1}v_1^{n-1}v_2^{m-1}([n,m]+[m,n])+\sum\limits_{n\geq 2, p+q=n-2}v_1^{p}v_2^{q}[n])$$
$$=\zeta (\sum\limits_{n,m\geq 1}v_1^{n-1}v_2^{m-1}[n]\qshuffle[m])$$
$$=\zeta\circ reg_V([v_1])\cdot \zeta\circ reg_V( [v_2]).$$

$$\zeta\circ reg_V([v_1,v_2]\shuffle [v_3])=\zeta\circ reg_V([v_1,v_2,v_3]+[v_1,v_3,v_2]+[v_3,v_1,v_2] )$$
$$=\zeta(\sum\limits_{n,m,p\geq 1}v_1^{n-1}v_2^{m-1}v_3^{p-1}[n,m]\shuffle [p]+$$
$$\sum\limits_{n\geq 2, p+q=n-2,m\geq 1}\left((\frac{1}{2}v_1^{p}v_2^{q}v_3^{m-1}+v_1^pv_2^{m-1}v_3^q)[n,m]+(\frac{1}{2}v_1^{p}v_2^{q}v_3^{m-1}+v_1^{m-1}v_2^pv_3^q)[m,n]
\right)+$$
$$\frac{1}{2}\sum\limits_{n\geq 3,\ p+q+r=n-3}v_1^{p}v_2^qv_3^r[n])$$
$$=\zeta(\sum\limits_{n,m\geq 1}v_1^{n-1}v_2^{m-1}[n,m]+
\frac{1}{2}\sum\limits_{n\geq 2, p+q=n-2}v_1^{p}v_2^{q}[n])\qshuffle 
\sum\limits_{r\geq 1}v_3^{r-1}[r])$$
$$=\zeta\circ reg_V([v_1,v_2])\cdot \zeta\circ reg_V([v_3]).$$

\end{document}